\documentclass[12pt,a4paper]{article}
 \usepackage{pdfpages}
\usepackage{parskip}
\usepackage{graphicx}
\usepackage[]{epsfig}
\usepackage{amsmath}
\usepackage{amssymb}
\usepackage{verbatim}
\usepackage{url}
\usepackage{latexsym}
\usepackage{amsmath}
\usepackage{color}
\numberwithin{equation}{subsection}
\numberwithin{equation}{section}
\begin{subequations}
\newtheorem{proposition}{Proposition}
\newtheorem{theorem}[subsection]{Theorem}

\end{subequations}
\date{}
\title{\bf Dynamics of a harvested prey-predator model with
 prey refuge depended\\ on both species}
\author{Md. Manarul Haque$^{\dagger}$ and Sahabuddin Sarwardi\footnote{Corresponding author}\\
Department of Mathematics, Aliah University \\
IIA/27, New Town, Kolkata - 700 156, India \\
e-mail: \textcolor{blue}{s.sarwardi@gmail.com} \\}
\begin{document}

\maketitle
\baselineskip 0.18in

\begin{abstract}
\noindent The present paper deals with a prey-predator model with prey refuge proportion to both species and independent harvesting of each species. Our study shows that using refuge as control, it can break the limit circle of the system and reach the required state of equilibrium level. It is established the optimal harvesting policy. The boundedness, feasibility of interior equilibria, bionomic equilibrium have been determined. The main observation is that coefficient of refuge plays an important role in regulating the dynamics of the present system. Moreover the variation of the coefficient of refuge change the system from stable to unstable and vice-versa. Some numerical illustration are given in order to support of our analytical and theoretical findings.
\end{abstract}



\section{Introduction}
The growing needs for more food  and more energy have led to increase exploitation of natural resources. The problem related to multi-species fisheries have been drawing  attention of researcher in the recent years (cf. Clark \cite{Clark1, Clark2, Clark3, Clark4}, Mesterton-Gibbons \cite{Mesterton-Gibbons}). The economic progress and ecological balance always have a conflicting interests. Therefore, concerning the conservation for the long term benefits of humanity, there is a wide range of interest in the use of bionomic modeling to gain greater insight in the scientific management of renewable resources like fisheries and forestry. An excellent introduction of optimal management of renewable resources has been presented by (cf. Clark \cite{Clark1}). The techniques and issues associated  with the dynamic economic models of natural resources exploitation are further developed by Clark (cf. Clark \cite{Clark2, Clark3, Clark5}). Harvesting is one of the most important issue to control extinction and minimization of exploitation of renewable resources. Renewable harvesting policy is indisputably one of the major and interesting problems from ecological and economical point of view. There is also a realistic phenomena for surviving himself is refuge. The use of spatial refuges by the prey is one of the more relevant behavioral traits that affect the dynamics of predator prey systems. By using refuges some fraction of prey population are partially protected against predators. The existence of refuges has significant influence on the coexistence of the predator-prey systems. It has a strong impact for stabilizing and destabilizing the dynamical nature of the system. The study of consequence of prey refuge on the dynamic of prey-predator interaction can be recognized as a major effect in applied mathematics and theoretical ecology (cf. Hassell and May \cite{Hassell and May}, Hassell \cite{Hassell}, Holling \cite{Holling}, Hoy \cite{Hoy}, Smith \cite{Smith}). Some of the empirical and theoretical works have investigated the effect of refuge and gave a decision that the refuge used by prey have a stabilizing effect on prey predator interaction and prey species can be prevented from extinction by using this policy (cf. Collings \cite {Collings}, Freedman \cite{Freedman}, Gonzalez-Olivares and Ramos-Jiliberto \cite{Gonzalez-Olivares and Ramos-Jiliberto}, Hochberg and Holt \cite{Hochberg and Holt}, Huang et al. \cite{Huang et al.}, Krivan \cite{Krivan}, McNair \cite{McNair}, Ruxton \cite{Ruxton}, Sih \cite{Sih}, Taylor \cite{Taylor} ).

Also, Brauer and Soudak (cf. Brauer and Soudak \cite{Brauer and Soudak1, Brauer and Soudak2, Brauer and Soudak3, Brauer and Soudak4}) studied a class of predator prey models under constant rate of harvesting and under constant percentage of harvesting of both species simultaneously. They have shown that how to classify the possibilities of the quantitative behavior of solution to locate the set of initial values in which the trajectories of the solution approaches to either an asymptotic stable equilibrium or an asymptotic stable limit cycle. Dai and Tang \cite{Dai and Tang} studied the following predator-prey model in which two ecological interacting species are harvested independently with constant rates of harvesting:

\begin{eqnarray}
  \frac{dx}{dt} &=& rx(1-\frac{x}{k})-a\phi(x)y-\mu,\nonumber\\
  \frac{dy}{dt} &=& y(-d+ca\phi(x))-h,\label{eqn1}
\end{eqnarray}
where $x(t)$, $y(t)$ represents prey and predator species respectively. $r$, $k$, $d$, $a$, $c$, $\phi(x)$ are the intrinsic growth rate, carrying capacity of prey, death rate of predator, maximum per capita consumption rate of predator, conversion rate and general predator response function on prey respectively; $\mu$, $h$ are constant harvesting rates. The most crucial element in an interacting population model is the ``functional response"�the expression
that describes the rate at which the number of prey is consumed by a predator.
They have shown the complicated dynamics of the system. Kar (cf. Kar \cite{Kar}) studied the prey-predator model with prey refuge and individuals are subject to proportional harvesting rates. As far as our knowledge goes there is no work in predator-prey model system with independent harvesting and prey refuge proportional to both the species have done even consideration of prey refuge proportional to both species is one step closure to reality. Keeping this in mind an attempt is made in the present investigation to study the effect of refuge as well harvesting on a Holling Type-II prey-predator model. In this paper we take the following model:

\begin{eqnarray}
  \frac{dx}{dt} &=& rx(1-\frac{x}{k})-\frac{p(1-my)xy}{1 +ax(1-my)}-q_{1}E_{1}x=F_1(x, y),\nonumber\\
  \frac{dy}{dt} &=& \frac{ep(1-my)xy}{1 +ax(1-my)}-dy-q_{2}E_{2}y=F_2(x, y),\label{eqn1}
\end{eqnarray}
where $x(t)$, $y(t)$ represents the prey and predator populations respectively at any time t. $r$, $k$, $p$, $m$, $d$, $q_{1}$, $q_{2}$, are all positive constants and have their biological meanings, accordingly $r$, represent the intrinsic growth rate of prey species, $k$ is the carrying capacity of the prey in absence of predator and harvesting, $p>0$ is the maximum per capita consumption rate of predator, $e$ $(0<e<1)$ is the efficiency by which predator converts consumed prey into new predator, $d>0$ is the death rate of predator. $E_{1}\geq 0$, $E_{2}\geq 0$, denotes the harvesting effort for the prey and predator respectively. $q_{1}E_{1}x$ and $q_{2}E_{2}y,$ represent the catch of the respective species, where $q_{1}$, $q_{2},$ are catchability coefficients of the prey and predator species respectively. The present model incorporates a refuge proportional to both the prey-predator determines. i.e., $mxy$ from the predator species, where $m\in [0,\, 1]$ is a constant. Incorporation of prey refuge leaves  $(1-my)x$ of the prey available to be hunted by the predator.

\section{Some preliminaries}\label{preli}
\subsection{Existence and positive invariance}
Letting $X = (x,y)^t$, $F : \mathbb{R}^2 \rightarrow \mathbb{R}^2$, $F = (F_1,F_2)^t$,
the system (\ref{eqn1}) can be rewritten as $\frac{d{X}}{dt}=F({X})$. Here
$F_i \in C^{\infty}(\mathbb{R})$,\, $i=1,2.$
Since the vector function $F$ is a smooth function of the variables $(x,y)$ in the positive quadrant $\Omega=\{(x,y);x>0,\,y>0\}\subset \mathbb{R}^{2}_{+}$, local existence and uniqueness of the solution set hold.
\subsection{Boundedness}
\begin{proposition}  All the solutions of the model (\ref{eqn1}) are bounded uniformly.\end{proposition}
{\bf Proof.} Let us consider a function $\xi=x+\frac{y}{e}.$
For any  $\zeta=  d+q_{2}E_{2} >0$,

\begin{eqnarray}\frac{d\xi}{dt}+\zeta\xi &=& rx(1-\frac{x}{k})-q_{1}E_{1}x-\frac{1}{e}(dy+q_{2}E_{2}y)+\zeta \bigl( x+\frac{y}{e}\bigr) \nonumber\\
                                    &\leq& \frac{k}{4r}(r+\zeta-q_{1}E_{1})^2=\kappa >0. \end{eqnarray}
 By applying the theory of differential inequality (cf. Brickhoff and Rota \cite{Brickhoff and Rota 1982}), we have the following inequality $0<\xi(x,y)<\frac{\kappa}{\zeta}(1-e^{-\zeta})+\xi(0)e^{-\zeta t}<\max{\{\frac{\kappa}{\zeta},\xi(0)\}}$. Therefore,  $\underset{t\rightarrow +\infty}\limsup\, \xi(t)\leq \frac{\kappa}{\zeta}$ with last bound independent of initial conditions.
Hence, all the solutions of the system (\ref{eqn1}) starting from $\mathbb{R}^2_{+}$ evolve with respect to time and remain in the compact region $\mathbb{R}_{xy}={\{{(x,y)\in \mathbb{R}^2_{+}}: \xi(x,y)\leq\frac{\kappa}{\zeta}+\epsilon \}}$  for any $\epsilon>0.$

\subsection{Persistence }
 Persistence of a predator-prey model system plays an important role in mathematical ecology since the criteria of persistence for ecological systems is a condition that ensuring the long-term survival of all the species. Here we have shown the persistence by using average Lyapunov function (cf. Gard and Hallam \cite{Gard and Hallam}), together with its boundedness.\\
Considering the average Lyapunov function
$P(x,\,y)=x^{\rho_{1}} y^{\rho_{2}},$
 where $\rho_{1}$ and $\rho_{2}$ are undetermined positive constant. Let us define the function $\phi$ as follows:
  \begin{eqnarray}\Phi(x, \,y)&=&\frac{\dot{P}(x, \,y)}{P(x, \,y)}=\rho_{1} \frac{\dot{x}}{x}+\rho_{2} \frac{\dot{y}}{y}\nonumber\\
 &=&\rho_{1}\bigl(r(1-\frac{x}{k})-\frac{p(1-my)y}{1 +ax(1-my)}-q_{1}E_{1}\bigr) +\rho_{2}\bigl(\frac{ep(1-my)x}{1 +ax(1-my)}-d-q_{2}E_{2}\bigr).\nonumber\end{eqnarray}

 Now, $\Phi(0,0)=\rho_{1} r-\rho_{1} q_{1}E_{1}-\rho_{2}d-\rho_{2}q_{2}E_{2} > 0 $, if $\rho_{1} r>\rho_{1} q_{1}E_{1}+\rho_{2}d+\rho_{2}q_{2}E_{2}.$
 \begin{eqnarray}\Phi(x_{1}, \,0)&=&\rho_{1}\bigl(r(1-\frac{x_{1}}{k})-q_{1}E_{1}\bigr)+\rho_{2}\bigl(\frac{epx_{1}}{1 +ax_{1}}-d-q_{2}E_{2}\bigr)\nonumber\\
 &=& \rho_{2}\bigl(\frac{epx_{1}}{1 +ax_{1}}-d-q_{2}E_{2}\bigr)>0, \, \hbox{if $p> \frac{d+q_{2}E_{2}}{e}\bigl(a+\frac{1}{k(1-\frac{q_{1}E_{1}}{r})}\bigr).$} \nonumber\end{eqnarray}
Hence, the solution of the system (\ref{eqn1}) is permanent, if the conditions $\rho_{1} r>\rho_{1} q_{1}E_{1}+\rho_{2}d+\rho_{2}q_{2}E_{2}$ and  $p> \frac{d+q_{2}E_{2}}{e}\bigl(a+\frac{1}{k(1-\frac{q_{1}E_{1}}{r})}\bigr)$ are satisfied.
\section{Analysis of equilibria}
The equilibria of the system $(\ref{eqn1})$ are
(i) the trivial equilibrium $E^{0}(0,\,0)$;\, (ii) the axial equilibrium $E^{1}(x_{1}, \,0)=\bigl(k(1-\frac{q_{1}E_{1}}{r}),\, 0 \bigr),$ \, and the interior equilibrium point
(iii) $E^*(x^*, \,y^*)$.

 We are interested on the interior equilibrium point $E^*$, where $x^*$ and $y^*$ are satisfying the following system of equations:
 \begin{eqnarray}
   r(1-\frac{x^*}{k})-\frac{p (1-my^*)y^*}{1 +ax^*(1-my^*)}-q_{1}E_{1}=0,
   \frac{ep (1-my^*)}{1 +ax^*(1-my^*)x^*}-d-q_{2}E_{2}=0.
   \end{eqnarray}\label{eqn2}
 One can easily found that the interior equilibrium $E^*$ is feasible if the conditions
   $ E_{1}<\frac{r}{q_{1}}(1-\frac{x^*}{k})$ and  $ E_{2}<\frac{1}{aq_{2}}(ep-ad)$ are satisfied.

 \section{Local stability}\label{local}
 We now investigate the local asymptotically stability (LAS) of $(\ref{eqn1})$ around the feasible equilibrium points.
 The Jacobian matrix at an arbitrary point $(x, \,y)$ is
  $$J=
 \left[\begin{array}{cc} r-2\,\frac{rx}{k}-q_{1}E_{1}-\frac{py(1-my)}{(1+ax(1-my))^2}&-\frac{px(ax(1-my)^2+1-2my}{(1+ax(1-my))^2})\\
\noalign{\medskip}\frac{epy(1-my)}{(1+ax(1-my))^2}& -d-q_{2}E_{2}+\frac{e px(ax(1-my)^2+1-2my}{(1+ax(1-my))^2}) \end{array}\right].$$
Thus, the Jacobian matrix of the system $(\ref{eqn1})$ around the trivial equilibrium point $E^{0}(0, \,0)$ is
$$J_0=
 \left[\begin{array}{cc} r-\,q_{1}E_{1}&0\\
\noalign{\medskip}0& -d-q_{2}E_{2} \end{array}\right].$$
The eigenvalues for the steady state $E^{0}(0, \,0)$ are $(r-\,q_{1}E_{1})$ and $-(d+q_{2}E_{2}).$
The eigenvalue $r-\,q_{1}E_{1}$ is positive or negative according to $E_{1}<\frac{r}{q_1}$ or  $E_{1}>\frac{r}{q_1},$ i.e., the equilibrium point $E^{0}(0, \,0)$ is locally asymptotically stable (cf. Fig 1(a)) if the BTP (Biotechnical Productivity ) of the prey species $(\frac{r}{q_1})$ is less than the effort $(E_1)$ for the prey species. If $BTP> E_{1}$, then equilibrium point $E^0 (0,0)$ will be a saddle for the system (\ref{eqn1}).

The Jacobian matrix at the point $E^{1}(x_{1}, \,0)$ of the system $(\ref{eqn1})$ is
$$ J_1=
 \left[\begin{array}{cc} r-2\,\frac{rx_1}{k}-q_{1}E_{1}&-\frac{px_1}{(1+ax_1)}\\
\noalign{\medskip}0& -d-q_{2}E_{2}+\frac{epx_1}{(1+ax_1)} \end{array}\right].$$
The eigenvalues of $J_1$ are $\lambda_{1}=r-2\,\frac{rx_1}{k}-q_{1}E_{1}$ and $\lambda_{2}=-d-q_{2}E_{2}+\frac{epx_1}{(1+ax_1)}.$
The eigenvalues are negative if the condition $\bigl(1-\frac{d+q_{2}E_{2}}{epk-ak(d+q_{2}E_{2})}\bigr)<\frac{q_1 E_1}{r}<1$ is satisfied and hence $E^{1}(x_{1}, \,0)$ is locally asymptotically stable (cf. Fig 1(b)).\\

The Jacobian matrix around the interior equilibrium point $E^*(x^*, \,y^*)$ of the system (\ref{eqn1}) is
$$J_{*}=
 \left[\begin{array}{cc} r-2\,\frac{rx^*}{k}-q_{1}E_{1}-\frac{py^*(1-my^*)}{(1+ax^*(1-my^*))^2}&-\frac{px^*(ax^*(1-my^*)^2+1-2my^*}{(1+ax^*(1-my^*))^2})\\
\noalign{\medskip}\frac{epy^*(1-my^*)}{(1+ax^*(1-my^*))^2}& -d-q_{2}E_{2}+\frac{e px^*(ax^*(1-my^*)^2+1-2my^*}{(1+ax^*(1-my^*))^2}) \end{array}\right].$$
Both the eigenvalues of $J_*$ will be negative if the following conditions hold:
\begin{eqnarray}
eax^*(1-my^*)^2(p+y^*-px^*)+ey^*&>&p y^*(my^*-p-mex^*)\nonumber\\&~&+rx^*(1+ax^*(1-my^*))^2,\label{LAS1}\label{LAS1}\\
py^*(1-my^*)&>&(1+ax^*(1-my^*))^2. \label{LAS2}\end{eqnarray}
Hence, under this parametric conditions (\ref{LAS1})-(\ref{LAS2}) the interior equilibrium point $E^*(x^*, \,y^*)$ will be locally asymptotically stable.
 \begin{figure}[!tbp]
  \centering
  \begin{minipage}[b]{0.8\textwidth}
    \includegraphics[width=\textwidth]{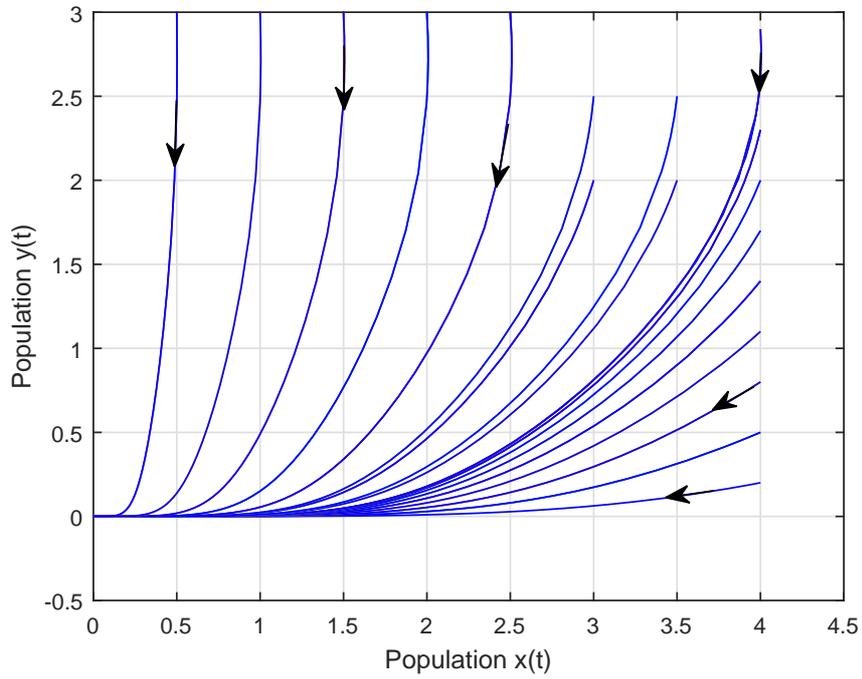}
    (a)
  \end{minipage}
  \hfill
  \begin{minipage}[b]{0.8\textwidth}
    \includegraphics[width=\textwidth]{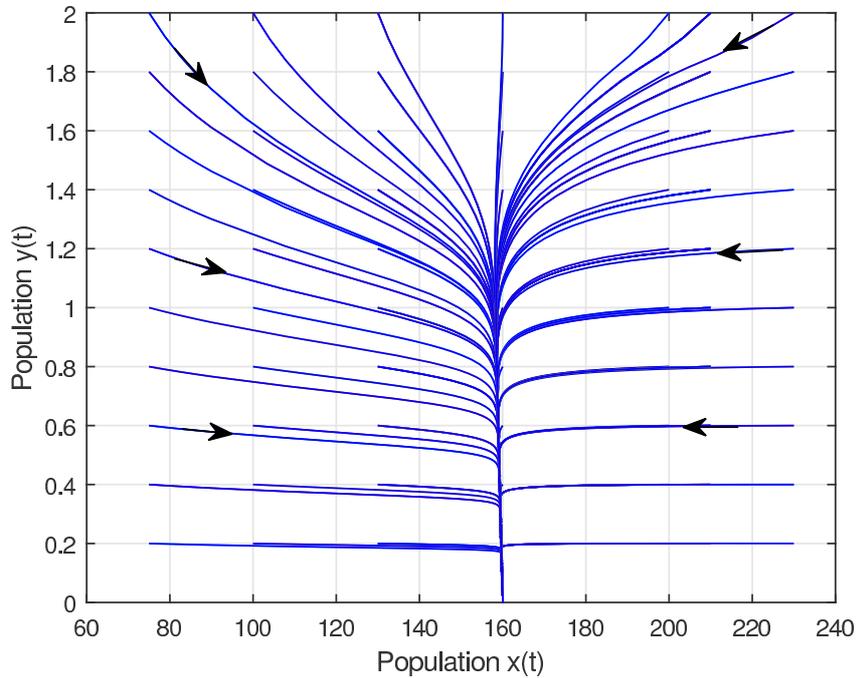}
    (b)
  \end{minipage}
 \caption{\small{\textrm (a) Demonstrates the local asymptotical stability (LAS) of the system (\ref{eqn1}) around the trivial equilibrium point $E^0(0, 0)$ corresponding to the parameter values: $r=1.0$, $a=0.04$, $d=0.5$, $m=0.5$, $p=0.2$, $q_1=0.4$, $q_2=0.6$, $E_1=3$, $E_2=1$, $k=200$, $e=0.25$. (b) Demonstrates the local asymptotical stability (LAS) of the system (\ref{eqn1}) around  the axial equilibrium point $E^1(x_1, 0)$ with the parameter values $r=3,$ $q_{1}=0.2$, remaining parameters are same as  (a).}}
\end{figure}

\section{Global stability around the interior equilibrium $E^{*}$}
\begin{theorem}The interior equilibrium $E^*$ is globally asymptotically stable if the condition $4rme(1+ak)(1+ax^*)>p m^2 y^*{^2}$ is satisfied.
\end{theorem}
{\bf Proof.}
Let us consider the suitable Lyapunov function as followers:\\
$V=(x-x^*-x^* ln\frac{x}{x^*})+(y-y^*-y^* ln\frac{y}{y^*})$.
Clearly V is positive definite for all $(x,\, y)\in \mathbf{R^2_{+}}\setminus (x^*, y^*)$.
Taking the time derivative along solution of the system (\ref{eqn1}),
we have
\begin{eqnarray}\frac{dV}{dt}&=&(x-x^*)\frac{\dot x}{x}+(y-y^*)\frac{\dot y}{y}\\
&=&-[\mathbf\alpha(x-x^*)^2+\mathbf\beta(y-y^*)^2+\mathbf\gamma(x-x^*)(y-y^*)],\label{eqg}
\end{eqnarray}
where $\alpha=\frac{r}{k}+\frac{pm(y+y^*)}{A}-\frac{py(1+myy^*)}{A},$  $\beta=\frac{mxep}{A},$ $\gamma=\frac{px(1+myy^*)+p(1-m(y+y^*))-ep-epy^*-p(y+y^*)}{A},$ $A=(1+ax(1-my))(1+ax^*(1-my^*)).$\\
Now, $4\alpha\beta-\gamma^2>0,$
if \begin{eqnarray} 4rme(1+ak)(1+ax^*)>p m^2 y^*{^2}\label{eqg1}.\end{eqnarray}\\
Thus, the  quadratic form (\ref{eqg}) is positive definite hence $\frac{dV}{dt}<0$ along all the trajectories in the first quadrant except $(x^*, \, y^*).$ Also $\frac{dV}{dt}|_{E^*}=0.$ The proof follows from the suitable Lyapunov function $V$ and Lyapunov-LaSalle's invariant principle (cf. Hale \cite{Hale}). Hence the equilibrium point is globally asymptotically stable if the condition (\ref{eqg1}) is satisfied.

\section{Bionomic equilibrium}
Let $c_1$ and $c_2$ are the fishing cost of prey and predator species per unit effort; $p_1$  and $p_2$ are the prices of prey and predator per unit biomass. \\
Therefore, the economic rent (net revenue) at any time can be taken as
$\pi = (p_1 q_1 x-c_1)E_1 +(p_2 q_2 x-c_2)E_2$
$ = \pi_x +\pi_y$ , where  $\pi _x=(p_1 q_1 x-c_1)E_1$ and $\pi_y=(p_2 q_2 y-c_2)E_2$
are the net revenue for the prey and predator species respectively.
The bionomic equilibrium is defined as a point where the biological and the economical equilibrium takes places.\\
The bionomic equilibrium $(x_{\infty}, \,y_{\infty}, \, E_1^{\infty}, \, E_2^{\infty})$ is given by the solutions of following simultaneous equations:
 \begin{eqnarray}
  r(1-\frac{x}{k})-\frac{p(1-my)y}{1 +ax(1-my)}-q_{1}E_{1}=0, \label{eqn3} \\
 \frac{ep(1-my)x}{1 +ax(1-my)}-d-q_{2}E_{2} =0, \label{eqn4} \\
 (p_1 q_1 x-c_1)E_1 +(p_2 q_2 y-c_2)E_2 =0. \label{eqn5}
 \end{eqnarray}

  The bionomic equilibria are determined  in different case as follows:\\

  $\bf Case I$: If $c_2 > p_2 q_2 y$, i.e., the fishing cost is greater than the revenue for the predator species, then the predator fishing will be stopped (i.e., $E_2$=0). Only the prey fishing will be in operational (i.e., $E_1< p_1q_1 x$ ).\\
   We then have $x_\infty=\frac{c_1}{p_1q_2}$. Therefore, putting this value in $(\ref{eqn3})$ and $(\ref{eqn4})$ we have the bionomic equilibrium point $(y_{\infty}, E_1^{\infty})$ in the y-$E_1$ plane if the parametric condition
   $E_1^{\infty} <\frac{r}{q_1} (1-\frac{c_1}{p_1q_1 k})$ holds.\\

   $\bf Case II$: If $c_1> p_1 q_1 x,$ i.e., the  fishing cost is greater than the revenue in the prey fishing then the prey fishery will be closed (i.e., $E_1=0$) . Only the predator fishing will be in operational (i.e., $c_2< p_2q_2 y$).
   In this case the bionomic equilibrium $(x_{\infty}, E_2^{\infty})$ will be in the first quadrant of x-$E_2$ plane if
   $E_2^{\infty}\geq \frac{1}{q_2} (d-\frac{p_2 q_2 e r k}{4 c_2}).$\\

   $\bf Case III$: If $c_1> p_1q_1 x$ and $c_2> p_2 q_2 y$, then the fishing cost is greater than the revenue for both the species and the whole fishery will be closed.\\

   $\bf Case IV$: If $c_1< p_1q_1 x$ and $c_2< p_2 q_2 y$ holds simultaneously, then the revenue for both the species  being positive then the whole fishery  will be in operational, in this case $x_{\infty} =\frac{c_{1}}{p_1 q_1}$, $y_{\infty} =\frac{c_2}{p_2 q_2}$ , therefore using this values in $(\ref{eqn3})$ and $(\ref{eqn4})$ we have the non trivial bionomic equilibrium point $(x_{\infty}, \,y_{\infty}, \, E_1^{\infty}, \, E_2^{\infty})$ exists if the parametric conditions
    (i) $\frac{r}{q_1} (1-\frac{c_1}{p_1 q_1 k})> \frac{p(p_2 q_2-m c_2)c_2}{p_2 q_2 \bigl(p_1 q_1 p_2 q_2 +a c_1(p_2 q_2 -m c_2)\bigr)}$\, and (ii) $d< \frac{ep (p_2 q_2 -m c_2)c_1}{p_1 q_1 p_2 q_2 +a c_1(p_2 q_2- mc_2)}$ are satisfied.

 \section{Optimal harvesting policy}
 The fundamental problem in commercial exploitation of renewable resources is to determine the optimal trade-off between current and future harvests. As observed by Clark (cf. Clark \cite{Clark1}). This problem which is very assured of resources conservation is an exceedingly different one not from the Mathematical view point perhaps, but certainly from a political and philosophical view point. However, we look at the problem from the economic view point only and we have to use the time discounting policy to handle the question of inter temporal benefits. This discounting is a normal practice in business management (cf. Solow \cite{Solow}). For determination of an optimal harvesting policy, we consider the present value $J^1$ of continuous stream of revenue as follows:
\begin{eqnarray}J^1=\int^\infty_{0} e^{-\delta t}[( p_{1}q_{1} x-c_1)E_1(t) + (p_{2}q_{2} y-c_2)E_2(t)]dt, \label{eqn6} \end{eqnarray}
where $\delta$ denote the instantaneous annual rate of discount. We are to optimize the equation $(\ref{eqn6})$ subject to the state equation $(\ref{eqn1})$ by using Pontryagian's maximum principle (cf. Pontryagin \cite{Pontryagin}).\\
Let us consider the Hamiltonian function $H$ as follows:\\
\begin{eqnarray} H &=& e^{-\delta t}\bigl((p_1 q_1 x-c_1)E_1+(p_2 q_2 y-c_2)E_2\bigr)+\lambda_{1}\bigl(rx(1-\frac{x}{k})-\frac{p(1-my)xy}{1 +ax(1-my)}-q_{1}E_{1}x\bigr) \nonumber\\ &~& +\lambda_{2}\bigl(\frac{ep(1-my)xy}{1 +ax(1-my)}-dy-q_{2}E_{2}y\bigr),\nonumber\end{eqnarray}
where $\lambda_1(t)$ and $\lambda_2(t)$ are adjoint variables. $E_1(t), E_2(t)$ are the control variables subject to the constraint $0\leq E_i(t)\leq (E_i)_ {max}$, $(i=1,2)$. The control variable $E_1(t), E_2(t)$  appear linearly in the Hamiltonian function $H$. Assuming that the control constraint are not binding, i.e., optimal solution does not occur at $(E_i)_ {max}$ or $(E_i)_ {min}$, we have singular control (cf. Clark \cite{Clark1}) given by $\frac{\partial H}{\partial E_i}=0.$\\
Therefore,
\begin{eqnarray} \frac{\partial H}{\partial E_1}=0, \implies  \lambda_1(t)= e^{-\delta t}(p_1-\frac{c_1}{q_1 x}); \,
\frac{\partial H}{\partial E_2}=0,\implies \lambda_2(t)= e^{-\delta t}(p_2-\frac{c_2}{q_2 y}). \label{eqn7} \end{eqnarray}\\
Thus, the shadow prices $e^{\delta t}\lambda_i(t)$, $(i=1,2)$ do not vary with time in optimal equilibrium. Hence they satisfy the transversally condition at $ t\rightarrow +\infty$, i.e., they remain bounded as $t \rightarrow \infty.$
Again, $\frac{\partial H}{\partial E_1}=0$ $\implies \lambda_1 q_1 x e^{\delta t}=\frac{\partial \pi_x}{\partial E_1}$;\, $\frac{\partial H}{\partial E_2}=0$ $\implies \lambda_2 q_2 y e^{\delta t}=\frac{\partial \pi_y}{\partial E_2}$.\\
From these relations one can say that for each species the user cost of harvesting per unit effort must be equal to the discounted value of the future marginal profit of effort at the steady state level.\\\\
Now we are to find out the optimal solution of the problem as follows:
\begin{eqnarray} \frac{d \lambda_1}{dt} &=&-\frac{\partial H}{\partial x} \nonumber\\
&=&-[e^{-\delta t}p_1 q_1 x E_1+\lambda_1{\{r(1-\frac{2x}{k})-\frac {py(1-my)}{(1+ax(1-my))^2}-q_1E_1\}}\nonumber\\&~&+
 \lambda_2 {\{\frac {epy(1-my)}{(1+ax(1-my))^2}\}}] \label{eqn8}
\end{eqnarray}
Substituting the values of $\lambda_1$ and $\lambda_2$ from $(\ref{eqn7})$, we have the relation between $x$ and $y$ as follows:
\begin{eqnarray} \delta \bigl(p_1-\frac{c_1}{q_1 x}\bigr)&=&p_1\Bigl(r\bigl(1-\frac{x}{k}\bigr)-\frac{py(1-my)}{1+ax(1-my)}\Bigr)+\bigl(p_2-\frac{c_2}{q_2 y}\bigr)\frac{epy(1-my)}{(1+ax(1-my))^2}\nonumber\\&~&+\bigl(p_1-\frac{c_1}{q_1 x}\bigr)\Bigl(-\frac{rx}{k}+\frac{paxy(1-my)^2}{(1+ax(1-my))^2}\Bigr)
.\label{eqn9}
\end{eqnarray}
Again \begin{eqnarray}
\frac{d \lambda_2}{dt} &=&-\frac{\partial H}{\partial y} \nonumber\\
&=& -\Bigl(e^{-\delta t} p_2 q_2 E_2+ \lambda_1 px\frac{(ax(1-my)^2+1-2my)}{(1+ax(1-my))^2}\Bigr)\nonumber\\
&~&+\lambda_2\Bigl(-d-q_2 E_2+p e x \frac{(ax(1-my)^2+1-2my)}{(1+ax(1-my))^2}\Bigr).\label{eqn10}
\end{eqnarray}
Substituting the values of $\lambda_1$ and $\lambda_2$ in $(\ref{eqn10})$ we have the following relation in $x$ and $y$
\begin{eqnarray}\delta (p_2-\frac{c_2}{q_2 y}) &=& p_2\Bigl(\frac{epx(1-my)}{1+ax(1-my)}-d\Bigr)+px(p_1-\frac{c_1}{q_1 x})\frac{(ax(1-my)^2+1-2my)}{(1+ax(1-my))^2} \nonumber\\&~&-(p_2-\frac{c_2}{q_2 y})\frac{epmxy}{(1+ax(1-my))^2}. \label{eqn11}
\end{eqnarray}
Solving equation $(\ref{eqn9})$ and $(\ref{eqn11})$ for $x$ and $y$ we have the optimal equilibrium $(x^*, y^*)$ and the optimal harvesting efforts $E_1^*$ and $E_2^*$ can be determined by the following equations:
\begin{eqnarray}E_1^{*} &=& \frac{1}{q_1}\bigl(r(1-\frac{x^*}{k})-\frac{(1-my^*)y^*}{1+ax^*(1-my^*)}\bigr),\nonumber\\
E_2^{*} &=& \frac{1}{q_2}\bigl(-d+ep\frac{(1-my^*)x^*}{1+ax^*(1-my^*)}\bigr).
\end{eqnarray}

 \section{Bifurcation analysis}
 The subject of bifurcation is the study of structurally unstable systems.
 This topic is a branch of mathematics concerned with dynamical systems which
 suffer sudden qualitative changes in parameters. A small change in parameter
 causes a topological change. The important question in the field of local bifurcation
 theory is that a system depends on a control parameters, as parameter changes, what
 happens to the non hyperbolic equilibria. Structurally unstable dynamic systems can be
 classified according to the number of parameters that appears in the differential
 equations describing the dynamics of the system. For good introduction on the basics
  of local bifurcation analysis the interested readers are referred to check out the books by  Guckenheimer \cite{Guckenheimer}, Wiggins \cite{Wiggins}, and Kuznetsov \cite {Kuznetsov}.
   In our study of bifurcation the parameter $r$ has been chosen for the transcritical bifurcation and $m$ has been taken as bifurcation parameter for the analysis of Hopf- bifurcation for the present system (\ref{eqn1}).

 \subsection{Existence of transcritical bifurcation around $E_1$}
The the system (\ref{eqn1}) experiences a transcritical
bifurcation around the axial equilibrium $E_1$ as the parameter
$r$ crosses its critical value $r=q_1 E_1=r_{tc}$. One of the eigenvalues
of $J_1$ will be zero iff $\text{Det}(J_1)=0$, which implies
either  $r-E_1q_1=0$, or $-d-q_2 E_2+\frac{epx_1}{1+ax_1}=0$. Let
$v$, $w$ are the eigenvector corresponding to the eigenvalue
$\lambda=0$ for the matrices $J_1$ and $J_1^{T}$. The eigenvectors
$v$, $w$ are found as $v=(v_1, v_2)^T=(1, 0)$ and $w=(w_1,
w_2)^T=\bigl(1, \,\frac{px_1}{epx_1-(1+ax_1)(d+E_2q_2)}\bigr)^T$. With these
eigenvectors, it is found that (i) $w^TF_r(E_1, \,
r_{tc})=0$,\,(ii)\,$ w^T[DF_r(E_1, \,
r_{tc})v]=(1-\frac{2x_1}{k})\neq 0,$\,(iii)\,$ w^T[D^2F(E_1, \,
r_{tc})(v, v)]=-\frac{2r_{tc}}{k}\neq 0.$
Hence, due to Satomayor theorem (cf. Sotomayar \cite{Sotomayor}),
the system experience transcritical bifurcation around the axial
equilibrium $E^1$ at $r=r_{tc}$.
 \subsection{Hopf bifurcation}
 It can be easily conclude that the equilibrium point $E^*$ may loss its stability through Hopf bifurcation under certain parametric condition. Considering $m$ as a bifurcation parameter one can detect the threshold value of $m = m_{h}$, which satisfy $\text{Tr}(J_{E^*})|_{m=m_{h}}=0$ and $\text{Det}(J_{E^*})|_{m=m_h}>0$.\\ The transversality condition for the Hopf bifurcation (cf. Carr \cite{Carr}, Hassard and Kazarinoff \cite{Hassard and Kazarinoff}, Perko \cite{Perko}) is $\frac{d}{dm}(\text{Tr}(J_{E^*}))|_{m=m_h}=\frac{py_m(2aex_m^2-amx_m y_{m}^2+ax_m y_m+2ex_m-y_m)}{(1+ax_m (-my_m+1))^3}|_{m=m_{h}}\neq0$, where $x_m$ and $y_m$ indicate their functionality with respect to the parameter $m$.\\
 The interior equilibrium point $E^*$ loss its stability through the non-degenerate Hopf-bifurcation when the parametric restriction $\text{Tr}(J_{E^*})|_{m=m_{h}}=0$ and the transversality conditions are satisfied simultaneously. \\
 Now we calculate the Lyapunov number to determine the nature of Hopf-bifurcating periodic solutions. Introducing perturbations $x=x_{1}+x_m|_{m=m_h},$ $y=y_{1}+y_m|_{m=m_h}$ in (\ref{eqn1}) and  then expanding in Taylor series, we have

 \begin{eqnarray}\dot{x_1}&=&a_{10} x_1+a_{01}y_1+a_{20}x_{1}^2+a_{11}x_1 y_1+a_{30}x_{1}^3+a_{21}x_{1}^2y_1+\cdots,\\
 \dot{y_1}&=&b_{10} x_1+b_{01}y_1+b_{20}x_{1}^2+b_{11}x_1 y_1+b_{30}x_{1}^3+b_{21}x_{1}^2y_1+\cdots,\end{eqnarray}

 where $a_{10}$, $a_{01}$, $b_{10},$ $b_{01}$ are the elements of the Jacobian matrix evaluated at the equilibrium point $E^*$ with $m=m_{h}$, hence $a_{10}+b_{01}=0$ and $\Delta=a_{10}b_{01}-a_{01}b_{10}>0.$
 The expression of the coefficients $a_{ij}$ and $b_{ij}$ are given bellow:\\

 $a_{10}=\frac{\partial F_1}{\partial x}|_{(E^*,\, m_h)}$,\, $a_{01}=\frac{\partial F_1}{\partial y}|_{(E^*,\, m_h)}$, $a_{12}=\frac{1}{2}\frac{\partial^3 F_1}{\partial x\partial y^2}|_{(E^*,\, m_h)}$,\, $a_{21}=\frac{1}{2}\frac{\partial^3 F_1}{\partial x^2\partial y}|_{(E^*,\, m_h)}$,\\

 $a_{20}=\frac{1}{2}\frac{\partial^2 F_1}{\partial x^2}|_{(E^*,\, m_h)}$,\, $a_{11}=\frac{\partial^2 F_1}{\partial x\partial y}|_{(E^*,\, m_h)}$, $a_{30}=\frac{1}{6}\frac{\partial^3 F_1}{\partial x^3}|_{(E^*,\, m_h)}$;\;

 $b_{10}=\frac{\partial F_2}{\partial x}|_{(E^*,\, m_h)}$,\, $b_{01}=\frac{\partial F_2}{\partial y}|_{(E^*,\, m_h)}$, $b_{12}=\frac{1}{2}\frac{\partial^3 F_2}{\partial x\partial y^2}|_{(E^*,\, m_h)}$,\, $b_{21}=\frac{1}{2}\frac{\partial^3 F_3}{\partial x^2\partial y}|_{(E^*,\, m_h)}$,\\

$b_{20}=\frac{1}{2}\frac{\partial^2 F_2}{\partial x^2}|_{(E^*,\, m_h)}$,\, $b_{11}=\frac{\partial^2 F_2}{\partial x \partial y}|_{(E^*,\, m_h)}$, $b_{30}=\frac{1}{6}\frac{\partial^3 F_2}{\partial x^3}|_{(E^*,\, m_h)}$.

 The value of first Lyapunov number (cf. Perko \cite{Perko}), which helps to determine the nature of the stability of limit cycle arising through Hopf-bifurcation is given by

 \begin{eqnarray}\sigma&=&-\frac{3\pi}{2 a_{10}\Delta^{\frac{3}{2}}}\Bigl[\Bigl(a_{10}b_{01}a_{11}^{2}+a_{10}a_{01}(b_{11}^2+a_{20}b_{11})-2a_{10}a_{01}a_{20}^2-a_{01}^2(2a_{20}b_{20}
 +b_{11}b_{20})\nonumber\\&~& -a_{11}a_{20} (a_{01}b_{10}-2a_{10}^2)\Bigr)-(a_{10}^2+a_{01}b_{10})\Bigl(-3a_{01}a_{30}+2a_{10}(a_{21}+b_{12})+(b_{10}a_{12}\nonumber\\&~&-a_{01}b_{21})\Bigr)\Bigr], \label{eqn12}\end{eqnarray} where the values of $a_{ij}$ and $b_{ij}$, $i,\,j = 0,1,2,3$ are included in {\bf Appendix A}.

 If the $\sigma<0$, the equilibrium point $E^*$ destabilized through a supercritical Hopf-bifurcation, and if $\sigma>0$ then the Hopf bifurcation is subcritical.

\section{Numerical simulations}
In this section we perform numerical simulation to validate our analytical findings of the previous Section by making uses of the computing Software Packages MATLAB-R2015a and Maple-18. The analytical findings of the present study are summarized and presented Schematically in Tables 1, Table 2 and Table 3. It is very difficult to validate the model results with realistic data so far proportional refuse and harvesting are considered in natural field. These results are verified by means of numerical illustrations of which some chosen ones are shown in figures. We took a hypothetical set of parameter values to illustrate our results. In Section $(7)$  optimal harvesting policy and corresponding effort are determined. For numerical justification we have taken a set of parameter values: $r=3$, $a=0.008,$ $d=0.04,$ $m=0.02,$ $p=0.2,$ $q_1=0.2,$  $q_2=0.6$, $k=500$, $e=0.15$, $p_1=2$, $p_2=3$, $c_1=1$, $c_2=2$, $\delta=.004$ for this set of parameter values the optimal equilibrium is found at $(188.5858, \, 30.6567)$ and the corresponding harvesting efforts are $E_1^{*}=1.8534$ and $E_2^{*}=5.8875.$ In Section $(8)$, it is found that the present system (\ref{eqn1}) experiences Hopf-bifurcation for the bifurcation parameter $m$ and we find out the threshold value of $m$ and Lyapunov number $\sigma$ to determine the nature of Hopf-bifurcation (cf. Sen et al. \cite{Sen}). For numerical validation we take the fixed set of parameter values: $r=3,$ $a=0.008,$ $d=0.04,$  $p=0.2$, $q_{1}=0.2$, $q_{2}=0.6,$  $E_{1}=2,$ $E_{1}=2,$ $k=500,$ $e=0.15$. It is found that the threshold value of the bifurcation parameter $m=m_{h}=0.010695$ at which the system (\ref{eqn1}) experiences Hopf bifurcation around the interior equilibrium $E^*$. For this choice of parameter values the first Lypunov number is $\sigma=-000143<0$. Hence the Hopf-bifurcation is supercritical. It is observed that the trivial and axial equilibria are locally asymptotically stable (LAS) by starting the solution plots from different initial conditions in the neighborhood of $E^{0}$ and $E^1$, all the solution plots eventually converge to the equilibria respectively (cf Fig. 1: (a)-(b)).  Fig. 2: (a)-(b) shows the limit cycle behavior of the system (\ref{eqn1}) around the equilibrium point $E^*=(67.86, \,18.00)$ for the set of parameter values: $r=3,$ $a=0.008,$ $d=0.04,$ $m=0.005$, $p=0.2$
$q_{1}=0.2$, $q_{2}=0.6,$  $E_{1}=2,$ $E_{1}=2,$  $k=500,$ $e=0.15$.
In Fig. 3 shows that the system experiences Hopf-bifurcation for bifurcation parameter $m=0.01$ and the other parameters are same as Fig. 2.
Fig. 4: (a)-(b)  shows that the system (\ref{eqn1}) is globally asymptotically stable  and conversing to the point $E^{*}(94.99, \, 23.33)$ for $m=0.015$ and the other parameters are same as Fig. 2. Fig. 5 Shows that gradual increase of the coefficient of refuge gives more protection to the prey species and reduce the rate of predation of predator, as a result the volume of prey species became larger and larger and the predator became smaller and smaller until it goes to extinct. This fact is presented
in tabular form (cf. Table 3).

  \begin{figure}
\centering
\begin{tabular}{cc}
(a) \epsfig{file=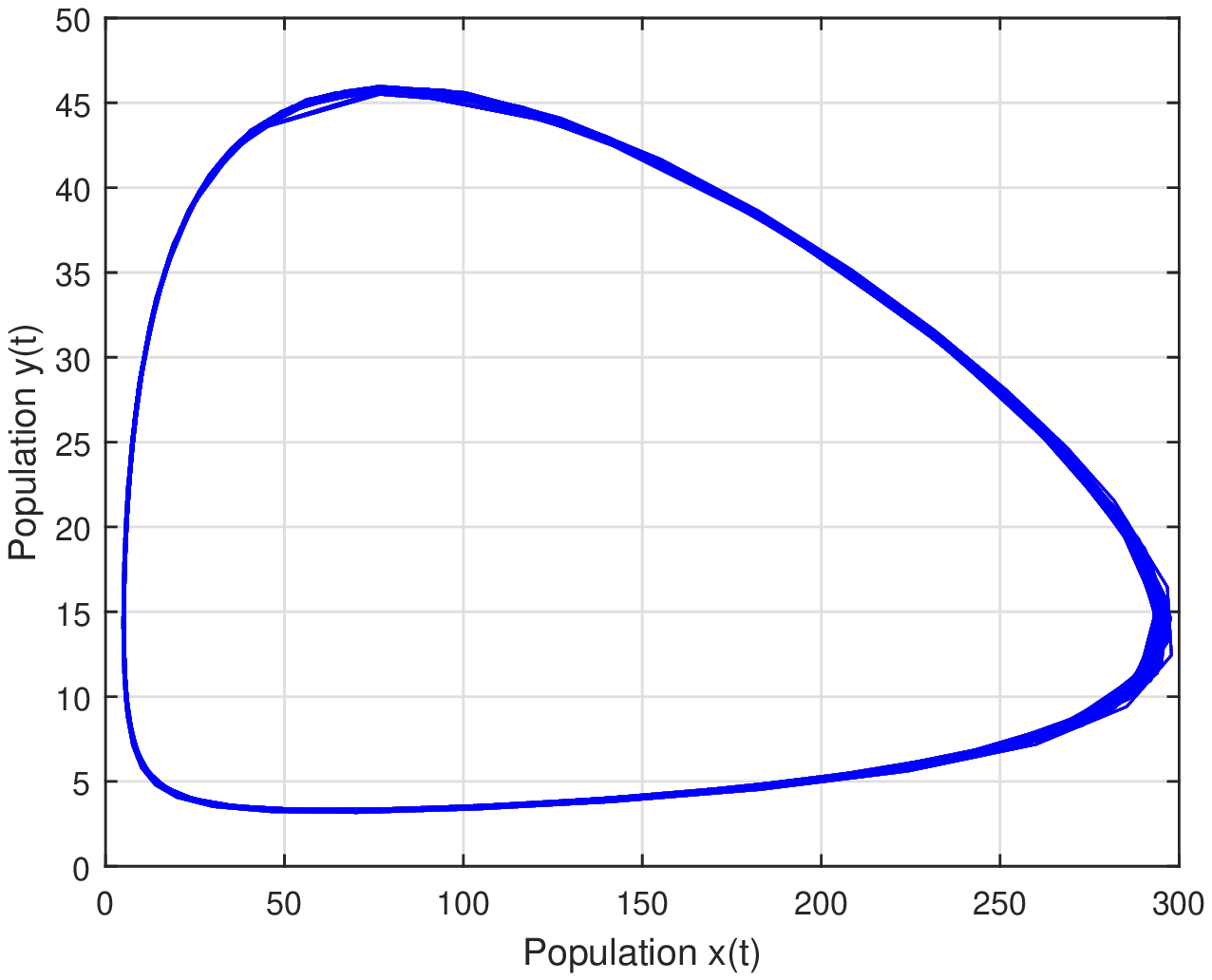,width=0.45\linewidth,height=6cm,clip=} &
(b) \epsfig{file=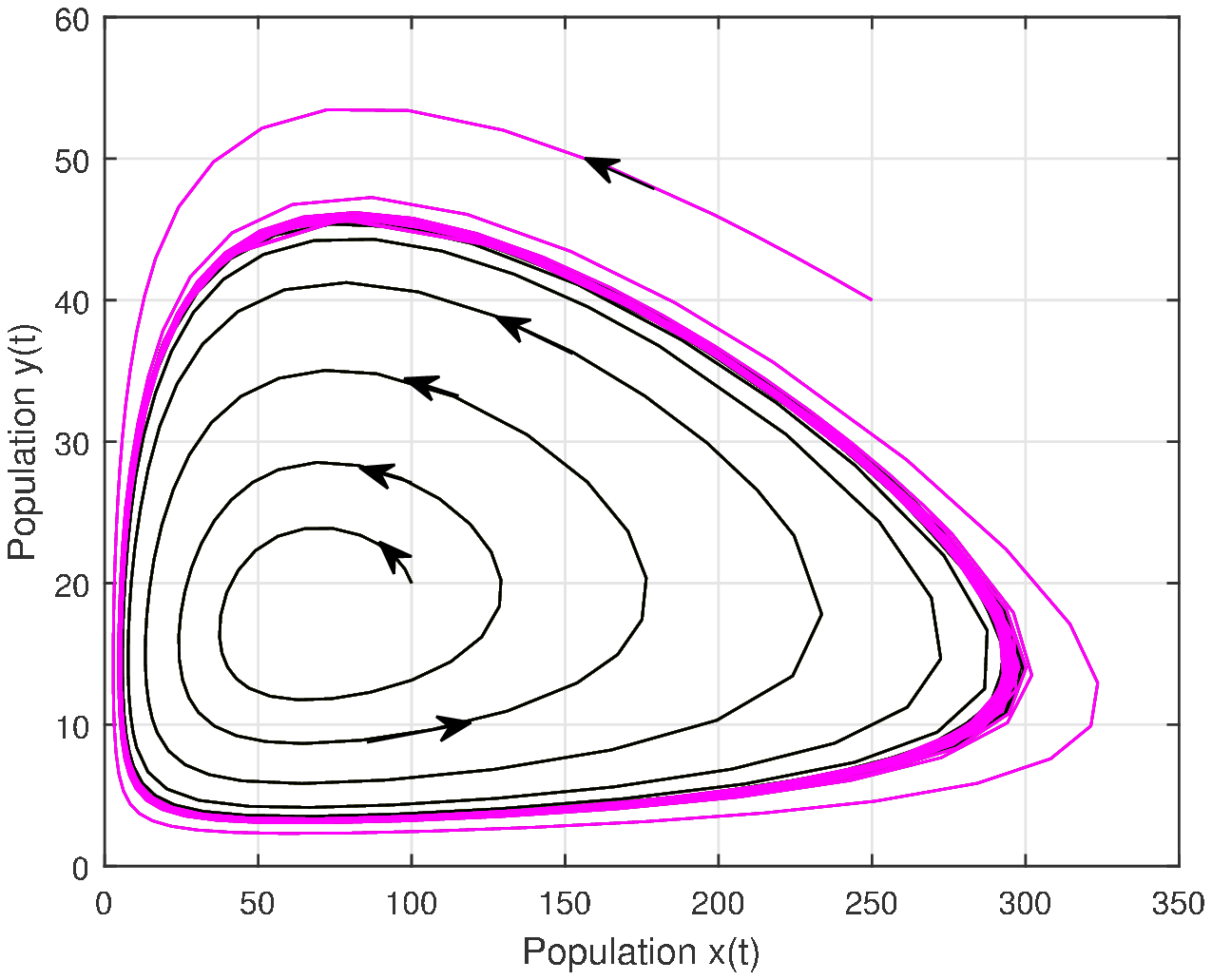,width=0.45\linewidth,height=6cm,clip=}
\end{tabular}
 \caption{\small{\textrm Solution plots $(a)$-$(b)$ depict that the system
(\ref{eqn1}) possesses a limit cycle surrounding  $E^{*} = (67.86, \, 18.00).$ Here the set of parameter
values used is:\, $r=3,$ $a=0.008,$ $d=0.04,$ $m=0.005$, $p=0.2$
$q_{1}=0.2$, $q_{2}=0.6,$  $E_{1}=2,$ $E_{2}=2,$
 $k=500,$ $e=0.15$.}}\label{f2}
\end{figure}

  \begin{figure}
\centering
\begin{tabular}{cc}
(a) \epsfig{file=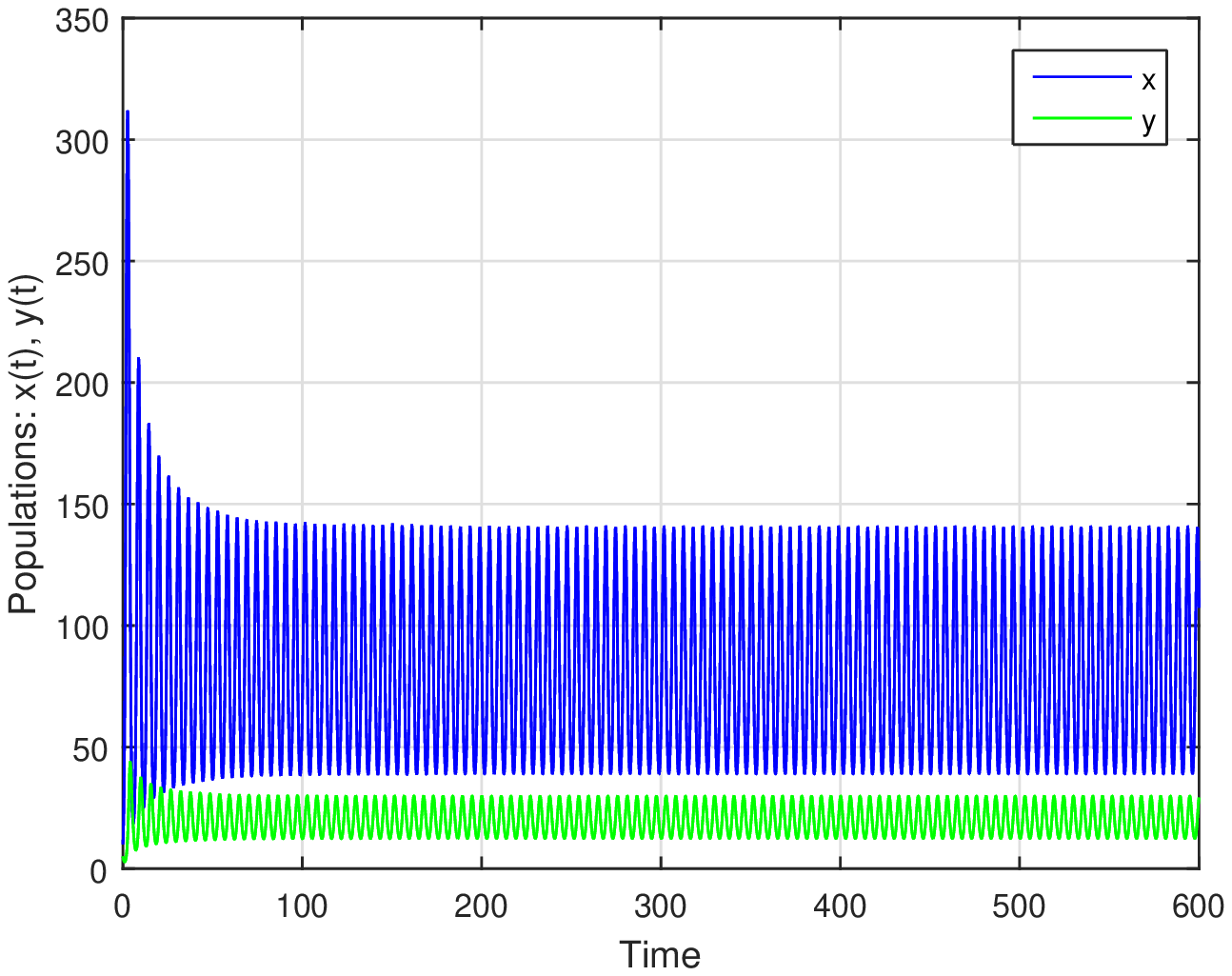,width=0.45\linewidth,height=6cm,clip=} &
(b) \epsfig{file=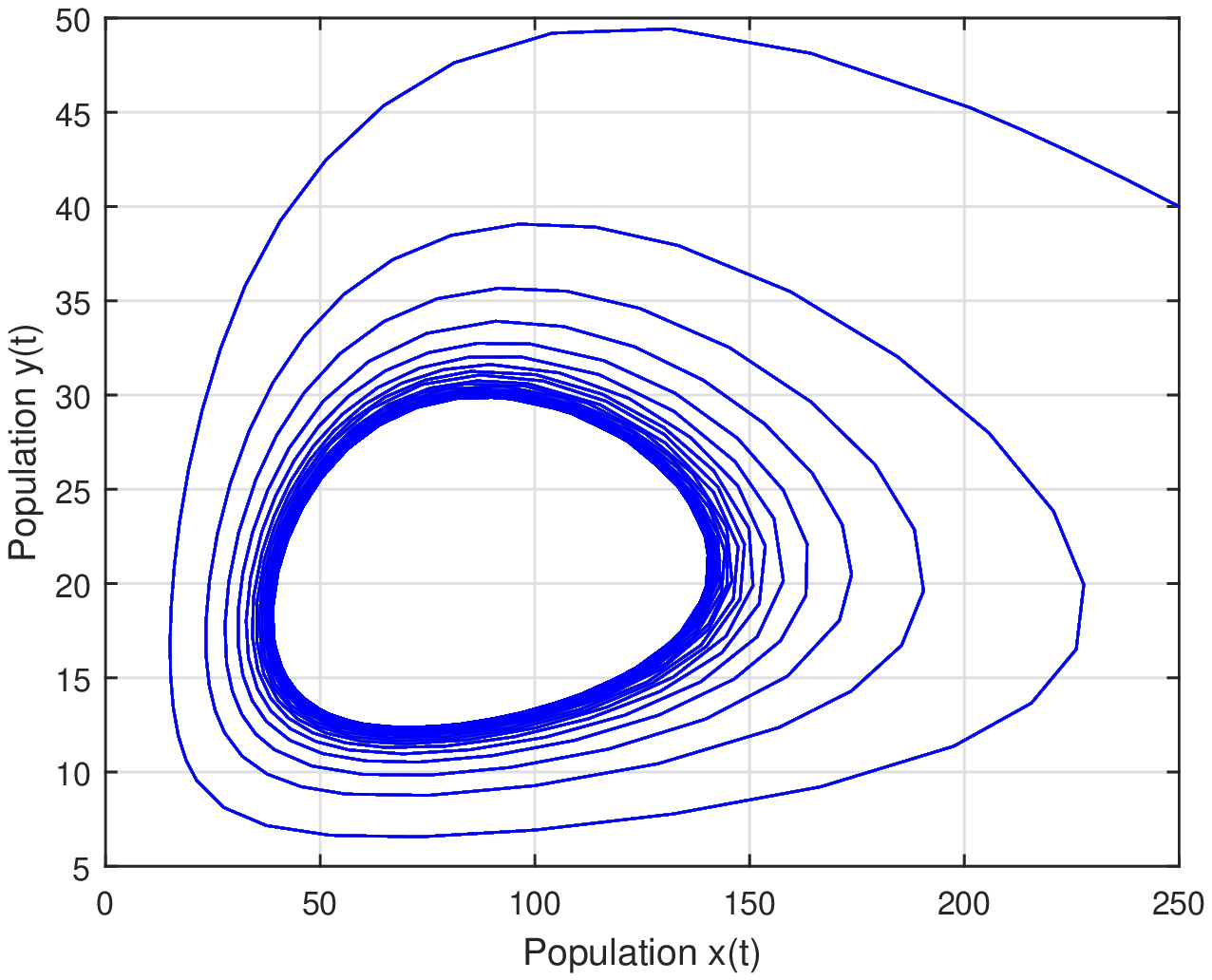,width=0.45\linewidth,height=6cm,clip=}
\end{tabular}
 \caption{\small{\textrm (a) There exists Hopf-bifurcating small amplitude periodic solutions. (b) the Phase diagram of the limit cycle for $m=0.01$ and the other parameter values used are same as Figure 2.}}\label{f3}
\end{figure}

\begin{figure}
\centering
\begin{tabular}{ccc}
(a) \epsfig{file=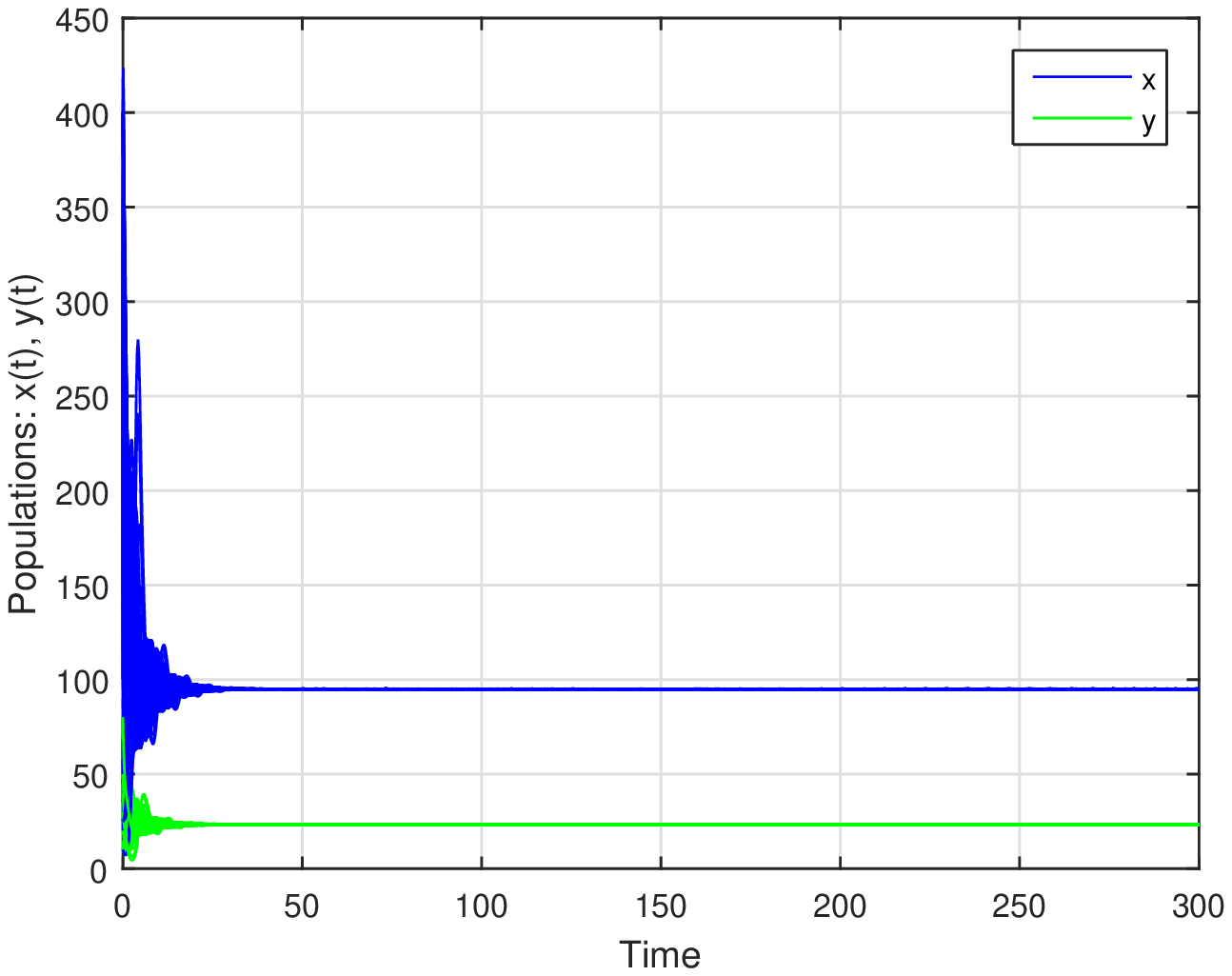,width=0.45\linewidth,height=5cm,clip=} &
(b) \epsfig{file=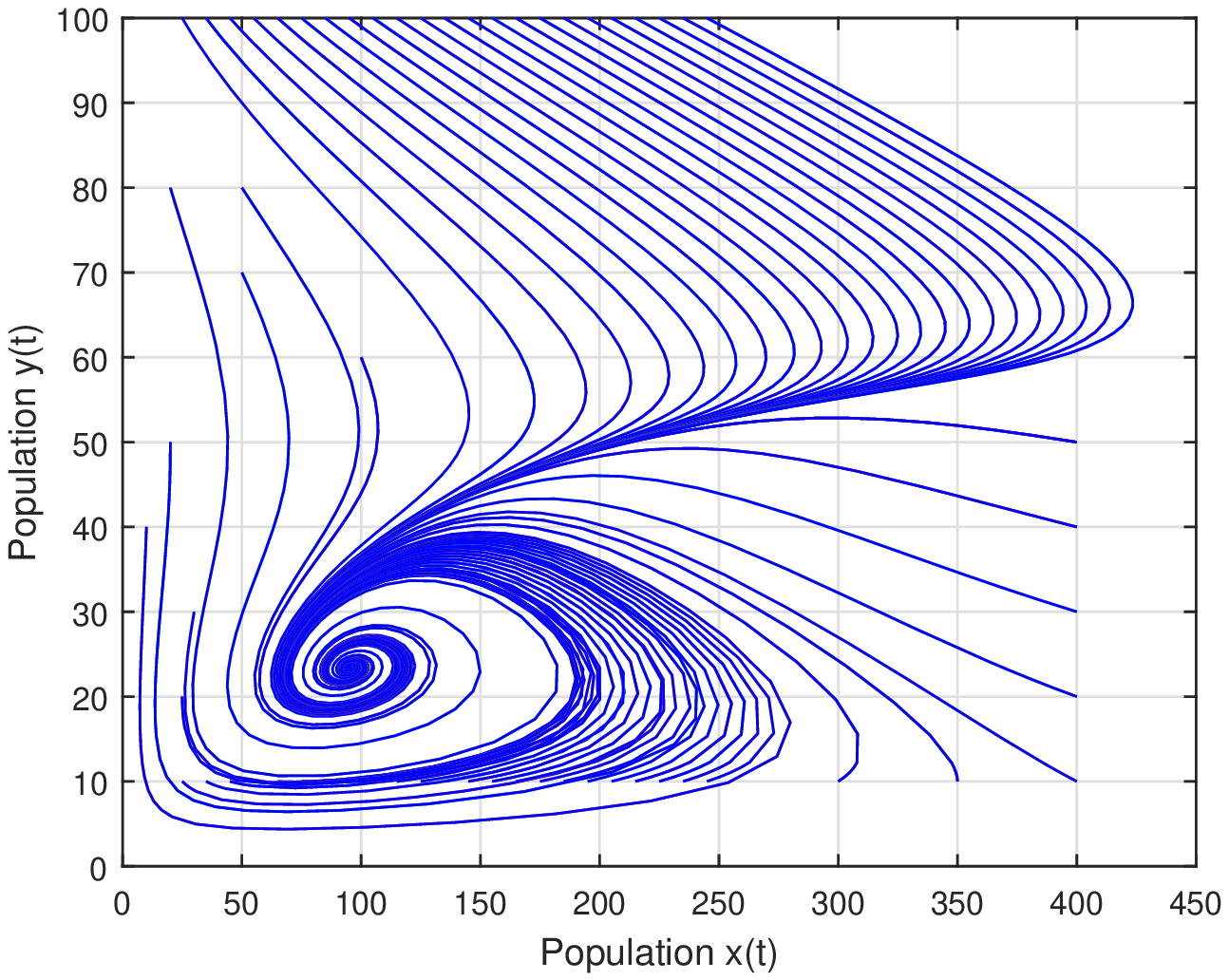,width=0.45\linewidth,height=5cm,clip=}&
\end{tabular}
\caption{\small{\textrm Solution plots of the system (\ref{eqn1}) converge to the interior equilibrium point $E^{*}(94.99, \, 23.33)$ for $m=0.015$. The other parameter values used are same as Figure 2.}}\label{f4}
\end{figure}

\begin{figure}[tbhp]\label{f5}
\begin{center}
\includegraphics[width=14cm, height=10cm]{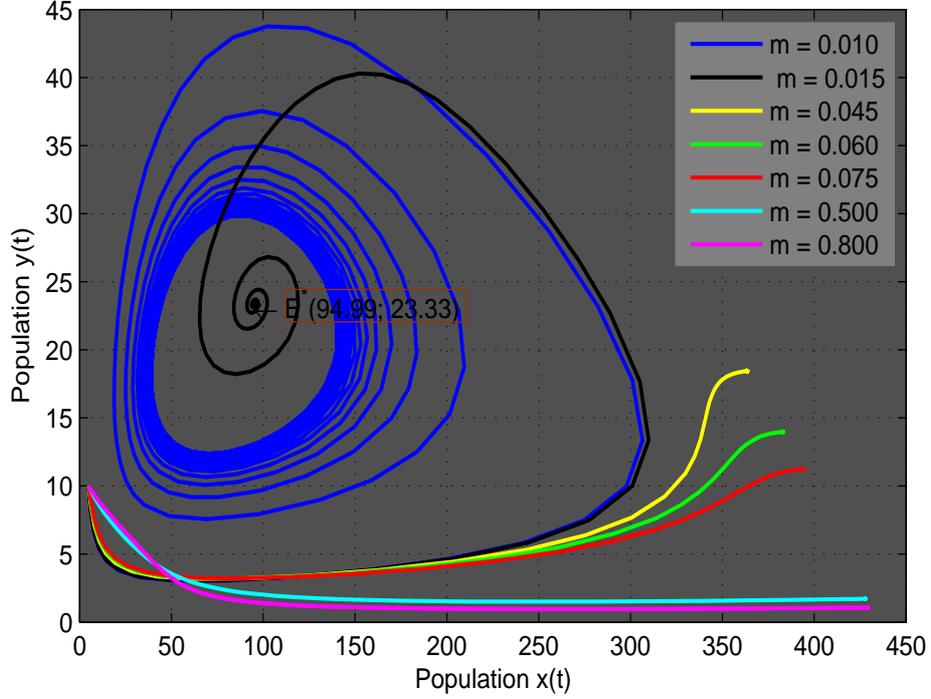}
\end{center}
\caption{\small{\textrm Different solution plots as $m$ increase (cf. Table \ref{Table3}) and other parameter values used are same as Figure 2.}}
\end{figure}

\begin{table}
\caption{\textrm {Schematic representation of the feasibility and stability conditions of the equilibria of the proposed model (\ref{eqn1}): LAS = Locally Asymptotically Stable; GAS = Globally Asymptotically Stable}}
\begin{center}
\begin{tabular}{|c|c|c|c|}\hline
{\textbf{Equilibria}}& \parbox[t]{1.5in}{\textbf{Feasibility conditions}}&\parbox[t]{1.80in}{\textbf{Stability/Persistent \mbox{}~~~conditions}}& \parbox[t]{1.1in}{\textbf{Nature}} \\
[.5ex]
\hline\hline
$E^{0}$& \parbox[t]{1.5in}{No condition} &\parbox[t]{1.80in} {$E_1<\frac{r}{q_1}$} & \parbox[t]{1.0in}{LAS}  \\
[.5ex] \hline
$E^{1}$ &\parbox[t]{1.5in}{$E_1<\frac{r}{q_1}$} &\parbox[t]{1.80in} {$\frac{r}{q_1}\bigl(1-\frac{d+q_2 E_2}{epk-ak(d+q_2E_2)}\bigr)<E_1<\frac{q_1}{r}$} &\parbox[t]{1.0in} {LAS}  \\
[.5ex] \hline

$E^{*}$&\parbox[t]{1.5in}{$E_1<\frac{r}{q_1}(1-\frac{x^{*}}{k})$,\, $E_2<\frac{1}{aq_2}(ep-ad)$} & \parbox[t]{1.80in} {\text{See the conditions} (\ref{LAS1})-(\ref{LAS2}) \text{of Section} \ref{local}} &\parbox[t]{1.0in}{LAS}  \\
[.5ex] \hline

$E^{*}$ &\parbox[t]{1.5in}{$E_1<\frac{r}{q_1}(1-\frac{x^{*}}{k})$,\, $E_2<\frac{1}{aq_2}(ep-ad)$} &\parbox[t]{1.80in} {$4rme(1+ak)(1+ax^*)>p m^2 y^*{^2}$} &\parbox[t]{1.0in} {GAS}  \\
[.5ex] \hline
\end{tabular}\label{Table1}
\end{center}
\end{table}

\begin{table}
\caption{Schematic representation of existence/feasibility conditions of the bionomic equilibria of the system (\ref{eqn1})}
\begin{center}
\begin{tabular}{|c|c|}
\hline
\parbox[t]{01.50in}{\quad\mbox{\bf Bionomic }\,\,\, {\bf Equilibrium}} & \parbox[t]{2.9in}{\bf{Existence/Feasibility conditions}}  \\
[1.5ex] \hline
$(y_{\infty}, E_1^{\infty})$ & \parbox[t]{2.9in}{$E_1^{\infty} <\frac{r}{q_1} (1-\frac{c_1}{p_1q_1 k}),$  provided prey fishing is in operational.}  \\
[1.6ex] \hline
$(x_{\infty}, E_2^{\infty})$ & \parbox[t]{2.9in}{$E_2^{\infty}\geq \frac{1}{q_2} (d-\frac{p_2 q_2 e r k}{4 c_2}),$  provided predator fishing is in operational.}  \\
[1.5ex] \hline
$(x_{\infty}, \,y_{\infty}, \, E_1^{\infty}, \, E_2^{\infty})$ & \parbox[t]{2.9in}{$\frac{r}{q_1} (1-\frac{c_1}{p_1 q_1 k})> \frac{p(p_2 q_2-m c_2)c_2}{p_2 q_2 \bigl(p_1 q_1 p_2 q_2 +a c_1(p_2 q_2 -m c_2)\bigr)},$\,\\ $d< \frac{ep (p_2 q_2 -m c_2)c_1}{p_1 q_1 p_2 q_2 +a c_1(p_2 q_2- mc_2)}$,\, \text{provided}\, both species fishing are in operational.}  \\
[1.5ex] \hline
\end{tabular}\label{Table2}
\end{center}
\end{table}

\begin{table}
\centering
\caption{The values of the components of equilibrium point $x^*$ and $y^*$ corresponding to the value of $m$.}
\begin{tabular}{ |c| c| c | c |c| c }\hline
No. & \parbox[t]{1.5in}{Fixed Parameters}& $m$ & \parbox[t]{0.9in}{$x^*,~y^*$}& {Figure 5} \\[3ex] \hline
 1 &\parbox[t]{2.1in}{$r_{1}=3, a=0.008, d=0.04, p=0.2,
q_{1}=0.2, q_{2}=0.6, E_1=2, E_2=2, k=500, e=0.15.$} &\parbox[t]{0.82in}{$0.010$} &\parbox[t]{1.3in}{$77.14,19.94$}& \parbox[t]{0.7in}{ blue}\\[2.5ex] \hline
 2&${\ldots\ldots\ldots\ldots\ldots\ldots\ldots}$ & \parbox[t]{0.82in}{$0.015$}& \parbox[t]{1.3in}{$94.99, 23.33$}& \parbox[t]{0.6in}{black} \\[1.5ex] \hline
  3& ${\ldots\ldots\ldots\ldots\ldots\ldots\ldots}$ & \parbox[t]{0.82in}{$0.045$} &\parbox[t]{1.3in}{$363.40, 18.45$} & \parbox[t]{0.6in}{yellow}  \\[0.5ex] \hline
  4 & ${\ldots\ldots\ldots\ldots\ldots\ldots\ldots}$ & \parbox[t]{0.82in}{$0.060$} &\parbox[t]{1.3in}{$383.05, 13.98$} & \parbox[t]{0.6in}{green}  \\[0.5ex] \hline
   5 & ${\ldots\ldots\ldots\ldots\ldots\ldots\ldots}$ & \parbox[t]{0.82in}{$0.075$} &\parbox[t]{1.3in}{$394.02, 11.24$} & \parbox[t]{0.6in}{red}  \\[0.5ex] \hline
    6 &${\ldots\ldots\ldots\ldots\ldots\ldots\ldots}$ &\parbox[t]{0.82in}{$0.500$} &\parbox[t]{1.3in}{$427.82, 1.71$}& \parbox[t]{0.6in}{cyan}\\[2.5ex] \hline
 7&${\ldots\ldots\ldots\ldots\ldots\ldots\ldots}$ & \parbox[t]{0.82in}{$0.800$}& \parbox[t]{1.3in}{$429.90, 1.07$}& \parbox[t]{0.6in}{magenta} \\[1.5ex] \hline

\end{tabular}\label{Table3}\\[3.5ex]
\end{table}

\section{Conclusion}
In this paper we consider a prey-predator harvesting model with Holling type -II response function incorporating prey refuge proportional to both the species. The novelty of our work lies in taking such kind of refuge function, which is more realistic phenomenon in ecosystem. We emphasize, on refuge coefficient that how the refuge function changes the system dynamics. We also observed from the both mathematical and empirical point of view that traits of behavioral policy of prey refuge has a stabilizing effect on a predator prey dynamics and this policy can help the prey species from extinction (cf. Anderson \cite{Anderson}, Cressman and Garay  \cite{Cressman and Gary}, Magalhaes et al. \cite{Magalhaes}, Rudolf et al. \cite{Rudolf}, Sarwardi et at. \cite{Sarwardi1}, Sarwardi et al. \cite{Sarwardi2}, Mukherjee \cite{Mukherjee}, Ma et al. \cite{Ma}). The key question is now the refuge alters the evolutionary dynamics of the system. In this paper we have shown that the dynamics of our system depends upon the functionality of refuge construction. In the present paper, it is assumed that functionality of prey refuge depends not only on prey size but on both the species. From the analytical view point the study of influence of prey refuge on the dynamics of interacting population is at present recognized as a significant and challenging issues (cf. Collings \cite{Collings}, Huang \cite{Huang}). We have shown that depending upon the bifurcation parameter $m$ the system exhibits stability as well as bifurcation around some of the equilibrium points. Also we have found that the optimal harvesting policy and the corresponding optimal effort using $E_1$  and $E_2$ as control parameter.

Before ending our conclusion we must say that, there are still some options to improve our model system to have much richer dynamics that what we have found in the present study. Here are given some rooms for our future studies. Firstly,
it would be more logical and proven reality that the term representing time delay used in digestion or gestation period for the predator species to produce new born have not been taken into account. Secondly, the harvesting efforts $E_{1} \text{and} E_{2}$ can be taken as time dependent functions. Thirdly, the harvesting effort can be taken as non-linear functions of both species. Fourthly, on the basis of the fact that more species could give more stable ecosystem, one can consider one more prey or predator or a pair of predator and prey into the exited system to have more stable system from the biological point of view.  The incorporation of all such relaxations existed in the present system in the future updated model would certainly be of some help to empirical researcher to predict their findings one step closer to the real situation from the ecological point of view. In our next paper we will study the dynamics of the system taking non-linear harvesting effort incorporating the taxation to sustain the renewable resources.

\skip0.5in
{\bf Acknowledgements:} Authors are thankful to the Department of Mathematics, Aliah University for providing opportunities to perform the present work. The corresponding author Dr. S. Sarwardi is thankful to his Ph.D. supervisor Prof. Prashanta Kumar Mandal, Department of Mathematics, Visva-Bharati (a Central University) for his generous help and continuous encouragement while preparing this manuscript.
\vskip0.5in

\section{Appendix}
 \begin{appendix}The expressions $a_{ij}$ and $b_{ij}$, appearing in (\ref{eqn12}) are defined by
 $a_{ij}=\frac{1}{i!j!}\frac{\partial ^{i+j}(F_{1})}{\partial x^i y^j}|_{(E^*,\, m_h)}$, \, $b_{ij}=\frac{1}{i!j!}\frac{\partial ^{i+j}(F_{2})}{\partial x^i y^j}|_{(E^*,\, m_h)}$ and their explicit values are as follows:\\
 $a_{10}=r-q_1E_1 -\frac{2x}{k}-\frac{ p(1-my) y}{(1+ax(1-my))^2}|_{E^*},$
 $a_{01}=-\frac{px\bigl(ax(1-my)^2+1-2my\bigr)}{(1+ax(-my+1))^2}|_{E^*},$
 $b_{10}=\frac{ep(1-my)y}{(1+ax(1-my))^2}|_{E^*}$,\\
 $b_{01}=-d-q_2E_2+\frac{epx\bigl(ax(1-my)^2+1-2my\bigr)}{(1+ax(1-my))^2}|_{E^*},$
 $a_{11}=\frac{-p(amxy-ax+2my-1)}{(1+ax(-my+1))^3}|_{E^*},$
 $a_{20}=-\frac{r}{k}+\frac{p(my-1)^2ya}{(1+ax(1-my))^3}|_{E^*}$,
 $a_{02}=\frac{-2pmx(ax+1)}{(1+ax(1-my))^3}|_{E^*},$
 $a_{30}=\frac{p(my-1)^3ya^2}{(1+ax(1-my))^4}|_{E^*},$ $b_{20}=\frac{ep(my-1)^2ya}{(1+ax(1-my))^3}|_{E^*},$\\
 $a_{21}=\frac{pa(my-1)(amxy-ax+3my-1)}{(1+ax(1-my))^4}|_{E^*},$
 $a_{12}= \frac{pm(a^2mx^2y-a^2x^2+2amxy+1)}{(1+ax(1-my))^4}|_{E^*},$$b_{30}=\frac{-ep(my-1)^3ya^2}{(1+ax(1-my))^4}|_{E^*},$
$b_{11}=\frac{ep(amxy-ax+2my-1)}{(1+ax(1-my))^3}|_{E^*},$
$b_{21}=\frac{-ep(my-1)a(amxy-ax+3my-1)}{(1+ax(1-my))^4}|_{E^{*}},$\\
$b_{12}=\frac{-epm(a^2mx^2y-a^2x^2+2amxy+1)}{(1+ax(1-my))^4}|_{E^{*}}$.
\end{appendix}

\end{document}